\documentclass{article}
\usepackage{amsfonts}

\title{Reducing Conjugacy in the full diffeomorphism 
group of $\R$ to conjugacy in the subgroup
of orientation-preserving maps}

\author{Anthony G. O'Farrell
\\
Mathematics Department\\
NUI, Maynooth\\
Co. Kildare\\
Ireland\\
\\
and
\\
\\
Maria Roginskaya
\\
Mathematics Department\\
Chalmers University of Technology and G\"oteborg University\\
SE-412 96 G\"oteborg\\
Sweden
}

\date{\today}

\begin{document}

\newtheorem{theorem}{Theorem}[section]%
\newtheorem{axiom}{Axiom}[section]%
\newtheorem{lemma}[theorem]{Lemma}%
\newtheorem{corollary}[theorem]{Corollary}%
\newtheorem{proposition}[theorem]{Proposition}%
\newtheorem{definition}{Definition}[theorem]
\newtheorem{definitions}[theorem]{Definitions}
\newtheorem{claim}{Claim}[theorem]
\newtheorem{example}{Example}[theorem]
\newtheorem{problem}{Problem}[section]

\def\Proof{{\par\medskip\noindent\bf Proof. }}
\def\qed{{\hfill\vrule height 4pt width 4pt depth 0pt
             \par\vskip\baselineskip}}

\def\R{{\mathbb{R}}}
\def\Z{{\mathbb{Z}}}
\def\N{{\mathbb{N}}}
\def\mod{{\,\hbox{mod}\,}}
\def\Diffeo{{\textup{Diffeo}}}
\def\Homeo{{\textup{Homeo}}}
\def\deg{{\textup{deg}}}
\def\fix{{\textup{fix}}}
\def\sign{{\textup{sign}}}
\def\Definition{{\par\medskip\noindent\bf Definition. }}
\def\ONE{1\hskip-4pt1}
\def\half{{\raise1pt\hbox{$\frac{\scriptscriptstyle1}{\scriptscriptstyle2}$}}}

\maketitle
\footnotetext{%
Mathematics Subject Classification: Primary 20E99, 
Secondary 20E36, 20F38, 20A05, 22E65, 57S25.  
\\Keywords: Diffeomorphism group, 
conjugacy, real line, orientation.%
\\Supported by SFI under grant RFP05/MAT0003.%
\\The authors are grateful to Ian Short for useful comments.}

\section*{Abstract}
Let $\Diffeo=\Diffeo(\R)$ denote the group of infinitely-differentiable
diffeomorphisms of the real line $\R$, under the
operation of composition, and let
$\Diffeo^+$ be the subgroup of diffeomorphisms
of degree $+1$, i.e. orientation-preserving diffeomorphisms.
We show how to reduce the problem of determining
whether or not two given elements $f,g\in \Diffeo$
are conjugate in $\Diffeo$ to associated conjugacy
problems in the subgroup $\Diffeo^+$.
The main result concerns the case when $f$ and $g$
have degree $-1$, and specifies (in an explicit and verifiable way)
precisely what must be added to the assumption
that their (compositional) squares are conjugate
in $\Diffeo^+$, in order to ensure that $f$ is
conjugated to $g$ by an element of $\Diffeo^+$.
The methods involve formal power series, and results of Kopell
on centralisers in the diffeomorphism group
of a half-open interval.

\section{Introduction and Notation}
Let $\Diffeo=\Diffeo(\R)$ denote the group of (infinitely-differentiable)
diffeomorphisms of the real line $\R$, under the
operation of composition.
In this paper we show how to 
reduce the conjugacy problem in $\Diffeo$
to the conjugacy problem in the index-two subgroup 
$$ \Diffeo^+=  \{f\in\Diffeo: \deg f = +1\},$$
where deg$f$ is the degree of 
$f$ ($=\pm1$, depending on whether or not
$f$ preserves the order on $\R$).

\medskip
We set some other notation:

$\Diffeo^-$:  $\{f\in\Diffeo: \deg f = -1\}$,
the other coset of $\Diffeo^+$ in $\Diffeo$.

$\Diffeo_0$:  the subgroup of $\Diffeo$ consisting
of those $f$ that fix $0$.

$\Diffeo^+_0$:  $\Diffeo_0\cap\Diffeo^+$.

fix$(f)$: the set of fixed points
of $f$.

$f^{\circ 2}$: $f\circ f$.

$f^{-1}$:
the compositional inverse of $f$.

$g^h$: $h^{-1}\circ g\circ h$, whenever
$g,h\in\Diffeo(I)$.  (We say that {\em $h$ conjugates $f$ to $g$}
if $f=g^h$.)

$-$: the map $x\mapsto -x$.

\medskip
We use similar notation for compositional powers and inverses
in the group $F$ 
of formally-invertible formal power series (with real coefficients)
in the indeterminate
$X$. The identity $X+0X^2+0X^3+\cdots$ is denoted simply by $X$.

$T_pf$ stands for the truncated Taylor series $f'(p)X+\cdots$ of a
function $f\in\Diffeo$. Note that $T_0$ is a homomorphism
from $\Diffeo_0$ to $F$, and 
$T_0(-) =-X$.

Typically, if $f$ and $g$ are conjugate
diffeomorphisms, then the family $\Phi$
of diffeomorphisms $\phi$ such that
$f=\phi^{-1}\circ g\circ\phi$
has more than one element. In fact
$\Phi$ is a left coset of the centraliser
$C_f$ of $f$ (and a right coset of $C_g$).
For this reason, it is important for us to
understand the structure of these centralisers.
The problem of describing $C_f$ is a special
conjugacy problem --- which maps conjugate
$f$ to itself?  Fortunately, this has already been 
addressed by Kopell [K].

\section{Preliminaries and 
Statement of Results}

\subsection{Reducing to conjugation by elements of $\Diffeo^+$}
The first (simple) 
proposition allows us to restrict attention to
conjugation using $h\in\Diffeo^+$.
\begin{proposition}\label{proposition-reduce-1}\label{general}
 Let $f,g\in\Diffeo$. Then the following
two conditions are equivalent:\\
(1) There exists $h\in\Diffeo$ such that $f=g^h$.
\\
(2) There exists $h\in\Diffeo^+$ such that $f=g^h$ or
$-\circ f\circ-=g^h$.
\end{proposition}
\Proof
If (1) holds, and $\deg h=-1$, then
$-\circ f\circ - = g^k$, with
$$ k(x) = h(-x).$$
The rest is obvious. 
\qed

\subsection{Reducing to conjugation of elements of $\Diffeo^+$}
The degree of a diffeomorphism is a conjugacy
invariant, so to complete the reduction of the
conjugacy problem in $\Diffeo$ to the problem
in $\Diffeo^+$, it suffices to deal with the
the case when $\deg f=\deg g=-1$ and
$\deg h=+1$.

\smallskip
{\it Let us agree that for the rest of this paper any objects
named $f$ and $g$ will be direction-reversing diffeomorphisms, and
any object named $h$ a direction-preserving diffeomorphism.}

\smallskip
Note that $\fix(f)$ and $\fix(g)$ are singletons.

If $f=g^h$,
then $h(\fix(f))=\fix(g)$,
and (since $\Diffeo^+$ acts transitively on $\R$)
we may thus, without loss in generality, suppose
that $f(0)=g(0)=h(0)=0$.

If $f=g^h$, then we also have
$f^{\circ 2} = 
(g^{\circ 2})^h$, 
$f^{-1}=(g^{-1})^h$,
and $f^{\circ 2}\in\Diffeo^+$.

We will prove the following reduction:
\begin{theorem}\label{theorem-reduction} Suppose $f,g \in
\Diffeo^-$, fixing $0$. Then the following two condition are
equivalent:
\begin{enumerate}
    \item $f=g^h$ for some $h\in \Diffeo^+$.
    \item \begin{enumerate}
        \item There exists $h_1\in \Diffeo^+_0$ such that $f^{\circ 2}
=(g^{\circ 2})^{h_1}$;

        and
        \item Letting $g_1=g^{h_1}$, there exists $h_2\in
        \Diffeo^+$, commuting with $f^{\circ 2}$ and fixing $0$,
        such that $T_0 f=(T_0 g_1)^{T_0 h_2}$.
    \end{enumerate}
\end{enumerate}
\end{theorem}

\subsection{Making the conditions explicit}
To complete the project of reducing
conjugation in $\Diffeo$ to conjugation in $\Diffeo^+$,
we have to find an effective way to check condition $2(b)$.
In other words, we have to replace the nonconstructive
\lq\lq there exists $h_2\in\Diffeo^+$"  
by some condition that can be checked algorithmically.
This is achieved by the following:

\begin{theorem}\label{theorem-last}
Suppose that $f,g\in\Diffeo^-$ both fix $0$,
and have $f^{\circ 2}
=g^{\circ 2}$. 
Then
there exists $h\in\Diffeo^+$, commuting with $f^{\circ 2}$,
such that $T_0f= (T_0g)^{T_0h}$
if and only if one of the following holds:
\begin{enumerate}
\item $(T_0f)^{\circ 2}\not=X$; 
\item $0$ is an interior point of $\fix(f^{\circ 2})$;
\item $(T_0f)^{\circ 2} = X$, $0$ is a boundary point
of $\fix(f^{\circ 2})$, and $T_0f=T_0g$.
\end{enumerate} 
\end{theorem}

Note that the conditions 1-3 are mutually-exclusive.
We record a couple of corollaries:

\begin{corollary}
Suppose $f,g \in
\Diffeo^-$, fixing $0$, and suppose $(T_0f)^{\circ2}\not=X$
or $0\in$int$\fix(f)$. Then 
    $f=g^h$ for some $h\in \Diffeo^+$
if and only if 
    $f^{\circ2}=(g^{\circ2})^h$ for some $h\in \Diffeo^+$.
\end{corollary}

In case $(T_0f)^{\circ2}\not=X$, any $h$ that conjugates $f^{\circ2}$
to $g^{\circ2}$ will also conjugate $f$ to $g$. In the other case
covered by this corollary, it is usually necessary to modify
$h$ near $0$.

\begin{corollary}
Suppose $f,g \in
\Diffeo^-$, fixing $0$, and suppose 
$(T_0f)^{\circ2}=X$
and $0\in$bdy$\fix(f)$. Then 
    $f=g^h$ for some $h\in \Diffeo^+$
if and only if 
    $f^{\circ2}=(g^{\circ2})^h$ for some $h\in \Diffeo^+$ and
$T_0f=T_0g$.
\end{corollary}

The last corollary covers the case where $0$ is
isolated in $\fix(f^{\circ2})$ and $T_0f$ is involutive,
as well as the case where $0$ is both an accumulation
point and a
boundary point of $\fix(f)$

\section{Proofs}

We begin by treating a special case:

\subsection{Involutions}
One possibility is that $f^{\circ 2}=\ONE$, i.e. $f$ is involutive,
and in that case so is any conjugate $g$. Conversely, we have:

\begin{proposition}\label{proposition-involution}
 If $\tau$ is a proper involution in $\Diffeo$, then it is
conjugated to $-$ by some $\psi\in \Diffeo^+$. Thus any two
involutions are conjugate.
\end{proposition}

\Proof
Let $\psi(x)=\half(x-\tau(x))$, whenever $x\in\R$.
It is straightforward to check that $\psi\in\Diffeo^+$, and
$\psi(\tau(x))=-\psi(x)$ for each $x\in\R$. Thus $\psi$
conjugates $\tau$ to $-$.
\qed

%
%

\subsection{Proof of Theorem \ref{theorem-reduction}}
\Proof.
$(1) \Rightarrow (2)$: Just take $h_1=h$ and $h_2=\ONE$.

\noindent
$(2)\Rightarrow (1)$: We just have to show that $f$ is conjugate to
$g_2=g_1^{h_2}$, and we note that 
$g_2^{\circ 2}=
(g_1^{\circ 2})^{h_2}
=f^{\circ 2}$.

Take
$$ k(x) = \left\{
\begin{array}{rcl}
x &,& x\ge0,\\
g_2(f^{-1}(x)) &,& x<0.
\end{array}
\right.
$$
Then, since $T_0f=T_0g_2$, we have $T_0(g_2\circ f^{-1})=X$,
so $k\in \Diffeo^+$.

We claim that $f=g_2^k$. Both sides are $0$ at $0$.

We consider the other two cases:

$1^\circ$, in which $x>0$. Then
    $$g_2^k(x)=k^{-1}(g_2(k(x)))=(g_2\circ f^{-1})^{-1}(g_2(x))=f(x).$$

$2^\circ$, in which $x<0$. Then
    $$g_2^k(x)=g_2(g_2(f^{-1}(x)))=
f^{\circ 2}(f^{-1})=f(x).$$

Thus the claim holds, and the theorem is proved. 
\qed

\subsection{The case when $f^{\circ 2}$ is not 
infinitesimally-involutive at $0$}
The nicest thing that can happen is that condition $(b)$ 
of Theorem \ref{theorem-reduction}
is automatically
true, once $(a)$ holds. 
The next theorem shows this does occur in a generic case 
(read $g_1$ for
$g$):

\begin{theorem}\label{theorem-8.4}
Suppose $f,g\in \Diffeo^-$, fixing $0$, with $f^{\circ 2}
=g^{\circ 2}$. 
Suppose  $(T_0f)^{\circ 2}\neq X$. Then $T_0f=T_0g$, and (by
Theorem \ref{theorem-reduction}) $f$ is conjugate to $g$.
\end{theorem}

Before giving the proof, we note a preliminary lemma:
\begin{lemma}\label{lemma-square}
The first nonzero term after $X$ in the (compositional) square
of a series with multiplier $-1$ has odd index.
\end{lemma}
\Proof
Let $S = -X +\cdots$ and $S^{\circ 2} = X$ mod $X^{2m}$.   
We claim that $S^{\circ 2} = X$ mod $X^{2m+1}$.  This will do.

Take $F = S-X$.  Then $F\circ S= S^{\circ 2} - S = -F$ mod $X^{2m}$,
so $F\circ  S\circ  F^{-1} = -X$ mod $X^{2m}$,
i.e. $F\circ  S\circ  F^{-1} = -X + cX^{2m}$ mod $X^{2m+1}$, for some $c\in\R$.
We calculate $F\circ  S^{\circ 2}\circ  F^{-1} = 
( F\circ  S\circ  F^{-1} )^{\circ 2} = 
X - c X^{2m} + c (-X)^{2m} = X$ 
mod $X^{2m+1}$,
so $S^{\circ 2} = X$ mod $X^{2m+1}$.
\qed

\medskip\noindent
{\bf Proof of Theorem \ref{theorem-8.4}.} 
\Proof
Let $q=g\circ f^{-1}$ and let $F=T_0f$, $G=T_0g$, and
$Q=G\circ F^{-1}=T_0(q)$. Then, since $$q^{-1}\circ
g=f=f^{\circ 2}f^{-1}=g\circ q$$ and $T_0$ is a group
homomorphism, we get
$$Q^{-1}\circ G=F=G\circ Q,$$
and deduce 
\begin{equation}\label{QFQ}
Q\circ F\circ Q=F
\end{equation}
and $F^{-1}\circ Q\circ F=Q^{-1}$, so that $Q$ is a
reversible series, reversed by $F$, 
and $Q$ commutes with $F^{\circ 2}$.

Note that (\ref{QFQ}) forces $Q=X (\mod X^2)$.

Now we consider the cases.

$1^\circ$. $f'(0)\neq -1$. Letting $\lambda=f'(0)$, there exists an
invertible series $W$ such that $F^W=\lambda X$. Letting $Q_1=Q^W$, we
see that $Q_1$ commutes with $\lambda^2 X$, and hence is $\mu X$ for
some nonzero real $\mu$. Since $Q_1=X (\mod X^2)$ also, we get
$\mu=1$, $Q_1=X$, $Q=X$, so $F=G$, and we are done.

$2^\circ$. $f'(0)=-1$. We may choose $p\in \N$ and a nonzero $a\in
\R$ such that
$$F^{\circ 2}=X+aX^{p+1} (\mod X^{p+2}).$$
Since $Q$ commutes with $F^{\circ 2}$, Lubin's Theorem [L, Cor.
5.3.2 (a) and Proposition 5.4] tells us that there is a $\mu\in \R$
such that $$Q=X+\mu X^{p+1} (\mod X^{p+2})$$ and if $\mu =0$ then
$Q=X$.

Suppose $\mu\not=0$. Then by Lemma \ref{lemma-square},
$p$ is even. But the first nonzero term after $X$
in a reversible series has even index (cf. [Ka], or [O, Theorem 5], 
for instance, or
calculate), so we have a contradiction. 
Hence, $\mu=0$, so $Q=X$, and we calculate again that
$F=G$, as in $1^\circ $.
\qed

\section{The case when $f^{\circ 2}$ is involutive on
a neighbourhood of $0$}




\begin{theorem}\label{theorem-8.6} Suppose $f,g \in \Diffeo^-$, fixing $0$, with
$f^{\circ 2}
=g^{\circ 2}$. 
Suppose $0$ is an 
interior point of $\fix (f^{\circ 2})$,
i.e. $f$ is involutive near $0$. Then there exists $h\in\Diffeo^+$,
commuting with $f^{\circ 2}$, fixing $0$, with $T_0f=(T_0g)^{T_0h}$,
and hence $f$ is conjugate to $g$.
\end{theorem}


\Proof Let $h_1(x)=\half(x-f(x))$, whenever $x\in\R$. Then
$h_1\in\Diffeo^+$, and $h_1(f(x)=-h_1(x)$ on $\fix(f^{\circ 2})$,
and hence on a neighbourhood of $0$.  Modifying $h_1$ off a 
neighbourhood of $0$, we may obtain $h_2\in\Diffeo^+$ 
with $h_2(x)=x$ off 
$\fix(f^{\circ 2})$.  It follows that $h_2$
commutes with $f^{\circ 2}$. 

Similarly, we may construct a function $h_3\in\Diffeo^+$
that commutes with $g^{\circ 2}=
f^{\circ 2}$ and has
$h_3(g(x))=-g(x)$ on a neighbourhood of $0$. Thus
$h= h_3^{-1}\circ h_2$ commutes with 
$f^{\circ 2}$
and has $h(f(x))=g(h(x))$ near $0$, so that
$T_0f=(T_0g)^{T_0h}$, as required.
\qed
\subsection{The Remaining Case}
We shall need the following
result from Kopell's paper [K, Lemma 1(b)]:
\begin{lemma}
Let $f,g\in\Diffeo^+$ both fix $0$ and commute.
If $T_0f=X$ and $0$ is not an interior point
of $\fix(f)$, then $T_0g=X$ as well.
\end{lemma}
\Proof\qed

\begin{theorem}\label{theorem-flat-special} Let $f,g \in\Diffeo^-$, fixing $0$, with
$f^{\circ 2}
=g^{\circ 2}$, and let $T_0f$ be involutive. Suppose 
$0$ is a boundary point of $\fix(f^{\circ 2})$.
Then $f$ is conjugate to $g$ if and only if
$T_0f=T_0g$.
\end{theorem}
\Proof
By Kopell's result,
any $h\in\Diffeo^+$ that commutes with
$f^{\circ 2}$ and fixes $0$ must have $T_0h=X$.
Thus the result follows from Theorem \ref{theorem-reduction}
\qed

Between them, Theorems \ref{theorem-8.4},
\ref{theorem-8.6} and
\ref{theorem-flat-special} cover all cases, and complete the proof of
Theorem \ref{theorem-last}.

\bibliographystyle{amsplain}

\section*{References}

\noindent
[L] J. Lubin. Nonarchimedean dynamical systems. Compositio
Mathematica 94 (1994) 321-46.

\noindent
[K] N. Kopell. Commuting diffeomorphisms. pp. 165-84 in
J. Palis + S. Smale (eds) Global Analysis. PSPM XIV. AMS. 1970.

\noindent
[Ka] E. Kasner. Conformal classification of analytic
arcs or elements: Poincar\'e's local problem of conformal geometry.
Transactions AMS 16 (1915) 333-49.

\noindent
[O] A.G. O'Farrell.  Composition of involutive power series, and
reversible series. Comput. Methods Funct. Theory 8 (2008) 173-93.

\bigskip
\noindent
e-mail:\\
anthonyg.ofarrell@gmail.com\\
maria@math.chalmers.se
\end{document}